\documentclass[10pt,reqno]{amsart}
\usepackage{lipsum}
\usepackage{amsmath}
\usepackage{mathtools}
\usepackage{enumerate}
\usepackage{amssymb}
\usepackage{tensor}
\usepackage{etoolbox}

\usepackage{tikz}
\usetikzlibrary{shapes.geometric}
\usetikzlibrary{cd}
\usetikzlibrary{calc}
\tikzcdset{every label/.append style = {font = \small}}

\usepackage{xr-hyper}
\usepackage{hyperref}

\newcommand{\B}{\mathcal{B}}
\newcommand{\C}{\mathcal{C}}

\newcommand{\I}{\mathcal{I}}

\newcommand{\T}{\textbf{T}}
\newcommand{\Tc}{\mathcal{T}}
\newcommand{\Sc}{\mathcal{S}}

\newcommand{\op}{^{\text{\normalfont op}}}

\newcommand{\e}{\epsilon}

\newcommand{\RTC}{\mathcal{RTC}}

\newcommand{\TC}{\mathcal{TC}}

\newcommand{\Ttwist}[1]{\Tr_{F,X}^{\mathcal{T}}(#1)}
\newcommand{\Stwist}[1]{\Tr_{F,X}^{\mathcal{S}}(#1)}

\DeclareMathOperator{\End}{End}

\DeclareMathOperator{\Hom}{Hom}

\DeclareMathOperator{\Vect}{\normalfont\underline{Vect}}
\DeclareMathOperator{\fld}{\mathbb{K}}
\DeclareMathOperator{\Yon}{\yen}

\DeclareMathOperator{\id}{id}
\DeclareMathOperator{\tr}{tr}
\DeclareMathOperator{\tid}{\normalfont \textbf{1}}

\DeclareMathOperator{\vprod}{\cdot}
\DeclareMathOperator{\Tr}{{\normalfont Tr}}

\newcommand{\tl}[1][\beta]{\text{\normalfont TL}}
\newcommand{\Rtl}[1][\beta]{\mathcal{R}\text{\normalfont TL}}
\newcommand{\tlred}[1][\beta]{\text{\normalfont TL}^\text{red}}
\newcommand{\Rtlred}[1][\beta]{\mathcal{R}\text{\normalfont TL}^\text{red}}

\newcommand{\Sharp}{^\sharp}

\newcommand{\defeq}{\vcentcolon=}

\newcommand{\SLn}[1]{\text{SL}_#1}

\newcommand\eo[2]{(#1#2,\epsilon_{#1}^{#2})}

\theoremstyle{plain}
\newtheorem{THM}{Theorem}[section]
\newtheorem{PROP}[THM]{Proposition}
\newtheorem{LEMMA}[THM]{Lemma}

\newtheorem{COR}[THM]{Corollary}
\theoremstyle{definition}

\newtheorem{DEF}[THM]{Definition}
\newtheorem*{DEF*}{Definition}

\title{Graphical Characterization of Modular Invariance}

\usepackage{amsaddr}
\author{Leonard Hardiman}

\email{hardiman@math.univ-lyon1.fr}

\begin{document}

\maketitle

\begin{abstract}

For a given modular tensor category we study representations of the corresponding tube category whose isomorphism classes are modular invariant matrices. In particular, we provide a characterization of these representations in terms of the annular graphical calculus of the tube category.

\end{abstract}

\section{Introduction}

Let $ \fld $ be an algebraically closed field, let $ \C $ be a modular tensor category over $ \fld $ together with a complete set of simples $ \I $ and let $ \tid \in \I $ denote the tensor identity. The \emph{modular data} of $ \C $ is composed of two $ \I \times \I $-matrices known as the S-matrix and the T-matrix respectively. These matrices satisfy the following equations,
	\begin{align*}
	(\Sc \Tc)^3 = \lambda \Sc^2, \ \Sc^2 = d(\C) C, \ \Tc C = C \Tc
	\end{align*}
where $ \lambda \in \fld $, $ d(\C) \in \fld $ is the \emph{dimension}~\cite[Definition 7.21.3]{Etingof15} of $ \C $ and $ C \defeq (\delta_{I^\vee,J})_{I,J \in \I} $ is the \emph{charge conjugation matrix}~\cite[Section 8.16]{Etingof15}. As $ C^2 = \id $ these equations imply that the modular data gives a projective representation of $ \SLn{2}(\mathbb{Z} $) a.k.a. the \emph{modular} group.

Part of the motivating force behind the study of modular tensor categories comes from their role in the mathematics of \emph{conformal field theory}. In particular, the category of modules over the vertex operator algebra of chiral symmetries naturally carries the structure of an modular tensor category. Let $ P $ be the partition function of the conformal field theory and let $ Z $ be the non-negative integer matrix describing the irreducible multiplicities of $ P $, i.e.\ $ P = \sum_{IJ} Z_{IJ}\chi_I \chi_J $ where the sum ranges over irreducible modules and $ \chi_I $ denotes the character of $ I $. Within this context the condition of conformal invariance for $ P $ may be rewritten as $ [Z,\Sc] = [Z,\Tc] = 0 $. This motivates the following definition.

	\begin{DEF} \label{def:modular_invariant}

	For a modular tensor category with tensor identity $ \tid $ and complete set of simples $ \I $, a \emph{modular invariant matrix} is a non-negative integer $ \I \times \I $-matrix that commutes with the modular data and whose $ (\tid,\tid) $-entry is $ 1 $\footnote{The condition that the $ (\tid,\tid) $-entry be $ 1 $ is also physically motivated, see~\cite[Section 5]{MR2164398}.}.

	\end{DEF}

The tube category, denoted $ \TC $, is a category whose objects coincide with those of $ \C $ and whose morphism are described by diagrams in $ \C $ inscribed upon the surface of a cylinder (see~\cite{hardiman_king} for a precise definition). Although $ \TC $ is not semisimple it may be though of as `semisimple in spirit' in the sense that its idempotent completion is semisimple. Indeed, for $ X $ and $ Y $ in $ \C $ we may consider the morphism $ \e_X^Y \in \End_{\TC}(IJ) $ given by~\ref{eq:idempotents}; in~\cite{hardiman2020graphical} it is shown that $ \e_X^Y $ is a idempotent and that the set $ \{\e_I^J\}_{I,J \in \I} $ forms a complete set of orthogonal, primitive idempotents in $ \TC $. The simple object corresponding to $ \e_I^J $ in the idempotent completion of $ \TC $ is denoted $ \eo{I}{J} $.

Let $ \Vect $ denote the category of finite dimensional vector spaces and let $ \RTC $ denote the category of contravariant functors from $ \TC $ to $ \Vect $. We call an object in $ \RTC $ a \emph{representation} of $ \TC $. As $ \TC $ admits a complete set of orthogonal, primitive idempotents the Yoneda embedding $ \Yon \colon \TC \to \RTC $ gives an alternative\footnote{Alternative, yet equivalent by \cite[Section 5.1.4]{Lurie09}.} idempotent completion of $ \TC $~\cite[Proposition 2.10]{hardiman2020graphical}. We use $ \eo{I}{J}\Sharp $ to denote the irreducible representation of $ \TC $ corresponding to the simple object $ \eo{I}{J} $, and, for $ F $ in $ \RTC $, we use $ F_I^J $ to denote the multiplicity space of $ \eo{I}{J}\Sharp $ in $ F $, i.e.
	\begin{align*}
	F_I^J \defeq \Hom_{\RTC}(\eo{I}{J}\Sharp,F).
	\end{align*}
Therefore, the isomorphism classes of objects in $ \RTC $ are in bijection with non-negative integer $ \I \times \I $-matrices. For $ F $ in $ \RTC $ we write $ Z(F) $ to denote the non-negative integer matrix corresponding to $ F $'s isomorphism class, i.e.\ $ Z(F)_{IJ} = \dim F_I^J $.

This naturally leads to the question of how to characterize objects in $ \RTC $ whose isomorphism class is a modular invariant matrix. We call such an object a \emph{modular invariant representation}. Furthermore, an object $ F $ in $ \RTC $ is called \emph{T-invariant} if it satisfies $ [Z(F),\Tc] = 0 $ and is called \emph{S-invariant} if it satisfies $ [Z(F),\Sc] = 0 $. The main result of this paper provides a characterization of modular invariant representations in terms of the graphical calculus of the tube category:

    \newtheorem*{cor:main_result}{Corollary \ref{cor:main_result}}
    \begin{cor:main_result}
    Let $ F $ be in $ \RTC $. Then $ F $ is modular invariant if and only if
        \begin{align}
         \Tr_{F,X}^{\Tc} = \Tr_{F,X} = \Tr_{F,X}^{\Sc} \label{eq:trace_equality}
        \end{align}
    for all $ X $ in $ \TC $.
    \end{cor:main_result}

Here $ \Tr_{F,X},\Tr_{F,X}^{\Tc},\Tr_{F,X}^{\Sc} \in \End_{TC}(X)^* $ are given by
$ \Tr_{F,X} \colon \alpha \mapsto \tr F(\alpha) $, \eqref{eq:t_twisted_trace} and \eqref{eq:s_twisted_trace} respectively. Alternatively, the values of $ \Tr_{F,X}^{\Tc} $ and $ \Tr_{F,X}^{\Sc} $ on $ \alpha \in \End(X) $ are given by taking $ \Tr_{F,X} $ and pre-composing it with an action of the standard generators of the mapping class group of the torus on the graphical description of $ \alpha $. Consequences of this interpreation are discussed in the concluding remarks found in Section~\ref{sec:conclusion}.

The structure of this paper is as follows: Section~\ref{sec:prilim} briefly introduces certain preliminaries on the graphical calculus of modular tensor categories and their associated tube categories. Section~\ref{sec:graph_char_t_inv} proves Corollary~\ref{cor:t_invariance_reformed} which states that the first equality in \eqref{eq:trace_equality} is equivalent to $ F $ being T-invariant. This section is fairly short as the corollary is a reformulation of a previous result from~\cite{hardiman2019extending}. Section~\ref{sec:graph_char_s_inv} proves Theorem~\ref{thm:s_invariance_reformulated} which states that the second equality in \eqref{eq:trace_equality} is equivalent to $ F $ being S-invariant. Finally, Section~\ref{sec:conclusion} combines these two results into Corollary~\ref{cor:main_result} and provides some concluding remarks.

\section{Preliminaries on Graphical Calculus} \label{sec:prilim}

For an overview of graphical calculus in modular tensor categories and their associated tube categories see~\cite[Sections 3-5]{hardiman2020graphical}. We now recall certain specific facts which will be used throughout this paper.

Firstly, we recall that the modular data of a modular tensor category $ \C $ admits the following graphical description,
   \begin{align} \label{eq:S_and_T_matrices}
   \Tc_{IJ} \defeq \delta_{I,J}
       \begin{array}{c}
           \begin{tikzpicture}[scale = 0.25]
           \draw[thick] (-1,1) -- (-1,3);
           \draw[thick] (-1,-1) -- (-1,-3);
           \draw[thick] (1,1) to[out = 90, in = 90] (3,1);
           \draw[thick] (1,-1) to[out = -90, in = -90] (3,-1);
           \draw[thick] (3,1) -- (3,-1);
           \draw[thick] (1,1) to[out=-90,in=90] (-1,-1);
           \draw[line width = 0.3cm, white] (-1, 1) to[out=-90,in=90] (1,-1);
           \draw[thick] (-1, 1) to[out=-90,in=90] (1,-1);
           \node at (-1, 4.3) {$ I $};
           \end{tikzpicture}
       \end{array}
   \quad
   \Sc_{IJ} \defeq
       \begin{array}{c}
           \begin{tikzpicture}[scale=0.5]
           \draw [thick](-1.5,0) arc (180:0:1);
           \draw[white,line width = 0.3cm] (0.5,0) circle (1);
           \draw[thick] (0.5,0) circle (1);
           \node at (-2,1) {$ I $};
           \node at (2,1) {$ J $};
           \draw [white,line width = 0.3cm](.5,0) arc (0:-180:1);
           \draw [thick](.5,0) arc (0:-180:1);
           \end{tikzpicture}
       \end{array}
   \end{align}
where we have exploited the fact that $ \End(I) = \fld $ for all $ I \in \I $. Secondly, for $ X $ and $ Y $ in $ \C $ we consider the following endomorphism,

	\begin{align} \label{eq:idempotents}
		\e_X^Y = \frac{1}{d(\C)} \bigoplus\limits_S d(S)
			\begin{array}{c}
				\begin{tikzpicture}[scale=0.15,every node/.style={inner sep=0,outer sep=-1}]
				\node (v1) at (0,5) {};
				\node (v4) at (0,-5) {};
				\node (v2) at (5,0) {};
				\node (v3) at (-5,0) {};
				\node (v5) at (-2.5,2.5) {};
				\node (v6) at (2.5,-2.5) {};
				\node (v7) at (1.5,3.5) {};
				\node (v11) at (3.5,1.5) {};
				\node (v12) at (-1.5,-3.5) {};
				\node (v10) at (-3.5,-1.5) {};
				\node (v13) at (2.5,-2.5) {};
				\draw [thick] (v7) edge (v10);
				\draw [line width =0.5em,white] (v5) edge (v13);
				\draw [thick] (v5) edge (v13);
				\draw [line width =0.5em,white] (v11) edge (v12);
				\draw [thick] (v11) edge (v12);
				\node at (2.75,4.75) {$ X $};
				\node at (4.75,2.75) {$ Y $};
				\node at (-3.75,3.75) {$ S $};
				\node at (3.75,-3.5) {$ S $};
				\draw[very thick, red]  (v1) edge (v3);
				\draw[very thick, red]  (v2) edge (v4);
				\draw[very thick]  (v1) edge (v2);
				\draw[very thick]  (v3) edge (v4);
				\end{tikzpicture}
			\end{array}
		\in \End_{\TC}(XY).
	\end{align}

Then (as mentioned in the introduction) $ \e_X^Y $ is a idempotent and that the set $ \{\e_I^J\}_{I,J \in \I} $ forms a complete set of orthogonal, primitive idempotents in $ \TC $. We use $ \eo{I}{J} $ to denote the corresponding simple object in the idempotent completion of $ \TC $.

\section{Graphical Characterization of $ \Tc $-invariance} \label{sec:graph_char_t_inv}
As before, let $ \fld $ be an algebraically closed field, let $ \C $ be a modular tensor category over $ \fld $ and let $ F $ be an object in $ \RTC $, i.e.\ $ F \colon \C\op \to \Vect $. We recall that $ F $ is called T-invariant if $ [Z(F),\Tc]=0 $. T-invariance was studied in~\cite[Section 4]{hardiman2019extending} where a graphical characterization was given in terms of the following automorphism in $ \TC $,
	\begin{align} \label{eq:twist}
	t_X \defeq
		\begin{array}{c}
			\begin{tikzpicture}[scale=0.15,every node/.style={inner sep=0,outer sep=-1}]
		 	\node (v1) at (0.5,4.5) {};
		 	\node (v2) at (4.5,0.5) {};
		 	\node (v3) at (-4.5,-0.5) {};
		 	\node (v4) at (-0.5,-4.5) {};
			\node (v6) at (-2,2) {};
			\node (v7) at (2.5,2.5) {};
			\node (v8) at (-2.5,-2.5) {};
			\node (v9) at (2,-2) {};
		 	\node at (4,4) {$ X $};
		 	\node at (-4,-4) {$ X $};
		 	\node at (3.5,-3.5) {$ X^\vee $};
		 	\node at (-3.5,3.5) {$ X^\vee $};
			\draw[thick] (v6) to[out=-45,in=-135] (v7);
			\draw[thick] (v9) to[out=135,in=45] (v8);
		 	\draw[very thick, red]  (v1) edge (v3);
		 	\draw[very thick, red]  (v2) edge (v4);
		 	\draw[very thick]  (v1) edge (v2);
		 	\draw[very thick]  (v3) edge (v4);
		 	\end{tikzpicture}
		\end{array}.
	\end{align}
	\begin{THM}[Theorem 4.4 in \cite{hardiman2019extending}] \label{thm:t_invariance}

	Let $ \C $ be a modular tensor category and let $ F $ be an object in $ \RTC $. $ F $ is T-invariant if and only if $ F(t_X) = \id_{F(X)} $ for all $ X $ in $ \C $.

	\end{THM}

The goal of this section is to rewrite the statement of this theorem into an equivalent statement. Although this reformulation may initially appear fairly arbitrary it illustrates the analogy with Theorem~\ref{thm:s_invariance}, our main result on S-invariance.

For $ X $ in $ \TC $ and $ F $ in $ \RTC $ we consider two elements in $ \End_{\TC}(X)^* $. One is simply $ \Tr_{F,X} \colon \alpha \mapsto \tr F(\alpha) $ and the other is given by
    \begin{align} \label{eq:t_twisted_trace}
    \Tr_{F,X}^{\Tc} \colon \alpha \mapsto \tr F(\alpha \circ t_X).
    \end{align}
We note that $ \Tr_{F,X}^{\Tc} $ does not depend on whether or not we compose $ t_X $ on the right or left by~\cite[Lemma 4.2]{hardiman2019extending}.

	\begin{COR} \label{cor:t_invariance_reformed}

	Let $ \C $ be a modular tensor category and let $ F $ be an object in $ \RTC $. $ F $ is T-invariant if and only if $ \Tr_{F,X}^{\Tc} = \Tr_{F,X} $ for all $ X $ in $ \C $.

	\proof

 	By~\cite[Proposition 4.3]{hardiman2019extending}, $ F(\e_I^J \circ t_{IJ}) = \frac{\Tc_{II}}{\Tc_{JJ}} F(\e_I^J) $. Therefore requiring that $ \tr F(\e_I^J \circ t_{IJ}) = \frac{\Tc_{II}}{\Tc_{JJ}} \tr F(\e_I^J) = \frac{\Tc_{II}}{\Tc_{JJ}} \dim F_I^J = \dim F_I^J = \tr F (\e_I^J) $ is equivalent to requiring that $ \dim F_I^J \neq 0 $ implies $ \frac{\Tc_{II}}{\Tc_{JJ}} = 1 $. As $ \Tc $ is diagonal, this condition (for all $ I,J \in \I $) is precisely the condition that $ Z(F) $ commutes with $ \Tc $, i.e.\ that $ F $ is T-invariant.

 	We therefore immediately have that $ \Tr_{F,X}^{\Tc} = \Tr_{F,X} $ implies T-invariance. We now suppose that $ F $ is T-invariant. Then, by Theorem~\ref{thm:t_invariance}, $ F(t_X) = \id_{F(X)} $ and we have $ \Ttwist{\alpha} = \tr F(\alpha \circ t_X) = \tr F(\alpha) = \Tr_{F,X} (\alpha). $

	\endproof

	\end{COR}

\section{Graphical Characterization of $ \Sc $-invariance} \label{sec:graph_char_s_inv}

The aim of this section is to provide analogues of Theorem~\ref{thm:t_invariance} and Corollary~\ref{cor:t_invariance_reformed} for S-invariance as opposed to T-invariance. We recall that $ F $ is called S-invariant if $ [Z(F),\Sc]=0 $. For $ I,J,S \in \I $ we consider the following endomorphism,
    \begin{align*}
	\gamma_{IJ}^S \defeq
		\begin{array}{c}
			\begin{tikzpicture}[scale=0.15,every node/.style={inner sep=0,outer sep=-1},rotate=-90]
			\node (v1) at (0,5) {};
			\node (v4) at (0,-5) {};
			\node (v2) at (5,0) {};
			\node (v3) at (-5,0) {};
			\node (v5) at (-2.5,2.5) {};
			\node (v6) at (2.5,-2.5) {};
			\node (v7) at (1.5,3.5) {};
			\node (v11) at (3.5,1.5) {};
			\node (v12) at (-1.5,-3.5) {};
			\node (v10) at (-3.5,-1.5) {};
			\node (v13) at (2.5,-2.5) {};
			\draw [thick] (v7) edge (v10);
			\draw [line width =0.5em,white] (v5) edge (v13);
			\draw [thick] (v5) edge (v13);
			\draw [line width =0.5em,white] (v11) edge (v12);
			\draw [thick] (v11) edge (v12);
			\node at (-5.25,-2.5) {$ I^\vee $};
			\node at (-2.75,-5.25) {$ J^\vee $};
			\node at (-4,4) {$ S $};
			\node at (4,-4) {$ S $};
			\draw[very thick]  (v1) edge (v3);
			\draw[very thick]  (v2) edge (v4);
			\draw[very thick, red]  (v1) edge (v2);
			\draw[very thick, red]  (v3) edge (v4);
			\end{tikzpicture}
		\end{array}
	\in \End_{\TC}(S)
    \end{align*}
    \begin{PROP} \label{prop:s_matrix_appears}

    For $ I,J,S \in \I $, we have
        \begin{align*}
        \gamma_{IJ}^S = \sum\limits_{A,B} \frac{\Sc_{I,A}\Sc_{B,J}}{d(A)d(B)} \e(S,A,B)
        \end{align*}
	where
		\begin{align*}
		\e(S,A,B) = \frac{d(A)d(B)}{d(\C)}\sum\limits_{T,b} d(T)
			\begin{array}{c}
				\begin{tikzpicture}[scale=0.2,every node/.style={inner sep=0,outer sep=-1}]
				\node (v1) at (0.5,5.5) {};
				\node (v4) at (-2,-7) {};
				\node (v2) at (5.5,0.5) {};
				\node (v3) at (-7,-2) {};
				\node (v5) at (-4,1) {};
				\node (v6) at (1,-4) {};
				\node (v7) at (0,2) {};
				\node (v11) at (2,0) {};
				\node (v12) at (-2,-4) {};
				\node (v10) at (-4,-2) {};
				\node (v13) at (1,-4) {};
				\draw [thick] (v7) edge (v10);
				\draw [line width =0.5em,white] (v5) edge (v13);
				\draw [thick] (v5) edge (v13);
				\draw [line width =0.5em,white] (v11) edge (v12);
				\draw [thick] (v11) edge (v12);
				\node at (-2.125,1.625) {$ A $};
				\node at (0,-0.5) {$ B $};
				\node at (-5,2) {$ T $};
				\node at (4,4) {$ S $};
				\node at (-5.5,-5.5) {$ S $};
				\draw[very thick, red]  (v1) edge (v3);
				\draw[very thick, red]  (v2) edge (v4);
				\draw[very thick]  (v1) edge (v2);
				\draw[very thick]  (v3) edge (v4);
				\node [draw,outer sep=0,inner sep=1,minimum size=10,minimum width = 26,fill=white,rotate=-45] (v9) at (1.5,1.5) {$ b $};
				\node [draw,outer sep=0,inner sep=1,minimum size=10,minimum width = 26,fill=white,rotate=-45] (v90) at (-3,-3) {$ b^* $};
				\node (v8) at (3,3) {};
				\draw [thick] (v8) edge (v9);
				\node (v14) at (-4.5,-4.5) {};
				\draw [thick] (v90) edge (v14);
				\end{tikzpicture}
			\end{array}
		\end{align*}
	where $ b $ ranges over a basis of $ \Hom_{C}(S,AB) $ and $ b^* $ is the dual element under the perfect pairing $ b \otimes b^* \mapsto \tr(b^* \circ b) $.
	By \cite[Corollary 5.7]{hardiman_king}, $ \e(S,A,B) $ is an idempotent corresponding (in the idempotent completion of $ \TC $) to the $ \eo{A}{B} $ isotypic component of $ S $.
    \proof As mentioned above, by \cite[Corollary 5.7]{hardiman_king}, we have $ \id_S = \sum_{A,B} \e(S,A,B). $
	Pre-composing this identify with $ \gamma_{IJ}^S $ yields
		\begin{align*}
		\gamma_{IJ}^S &=  \sum_{\substack{A,B,T \\ b}}\frac{d(A)d(B)d(T)}{d(\C)}
			\begin{array}{c}
				\begin{tikzpicture}[scale=0.2,every node/.style={inner sep=0,outer sep=-1}]
				\node (v1) at (0.5,5.5) {};
				\node (v4) at (-4.5,-9.5) {};
				\node (v2) at (5.5,0.5) {};
				\node (v3) at (-9.5,-4.5) {};
				\node (v5) at (-4,1) {};
				\node (v6) at (1,-4) {};
				\node (v7) at (0,2) {};
				\node (v11) at (2,0) {};
				\node (v12) at (-2,-4) {};
				\node (v10) at (-4,-2) {};
				\node (v13) at (1,-4) {};
				\node (v15) at (-3.5,-3.5) {};
				\node (v16) at (-7,-2) {};
				\node (v17) at (-2,-7) {};
				\node (v18) at (-8.5,-3.5) {};
				\node (v19) at (-3.5,-8.5) {};
				\node (v14) at (-7,-7) {};
				\draw [thick] (v16) edge (v17);
				\draw [line width =0.5em,white] (v15) edge (v14);
				\draw [thick] (v15) edge (v14);
				\draw [line width =0.5em,white] (v18) edge (v19);
				\draw [thick] (v18) edge (v19);
				\draw [thick] (v7) edge (v10);
				\draw [line width =0.5em,white] (v5) edge (v13);
				\draw [thick] (v5) edge (v13);
				\draw [line width =0.5em,white] (v11) edge (v12);
				\draw [thick] (v11) edge (v12);
				\node at (-2.125,1.625) {$ A $};
				\node at (0,-0.5) {$ B $};
				\node at (-5,2) {$ T $};
				\node at (-7.5,-0.75) {$ I^\vee $};
				\node at (-9.5,-2.5) {$ J^\vee $};
				\node at (4,4) {$ S $};
				\node at (-8,-8) {$ S $};
				\draw[very thick, red]  (v1) edge (v3);
				\draw[very thick, red]  (v2) edge (v4);
				\draw[very thick]  (v1) edge (v2);
				\draw[very thick]  (v3) edge (v4);
				\node [draw,outer sep=0,inner sep=1,minimum size=10,minimum width = 26,fill=white,rotate=-45] (v9) at (1.5,1.5) {$ b $};
				\node [draw,outer sep=0,inner sep=1,minimum size=10,minimum width = 26,fill=white,rotate=-45] (v90) at (-3,-3) {$ b^* $};
				\node (v8) at (3,3) {};
				\draw [thick] (v8) edge (v9);
				\end{tikzpicture}
			\end{array}\\
		&= \sum_{\substack{A,B,T \\ b}}\frac{d(A)d(B)d(T)}{d(\C)}
			\begin{array}{c}
				\begin{tikzpicture}[scale=0.2,every node/.style={inner sep=0,outer sep=-1}]
				\node (v1) at (0.5,5.5) {};
				\node (v4) at (-5,-10) {};
				\node (v2) at (5.5,0.5) {};
				\node (v3) at (-10,-5) {};
				\node (v5) at (-3.75,1.25) {};
				\node (v7) at (0,2) {};
				\node (v11) at (2,0) {};
				\node (v12) at (-5,-7) {};
				\node (v10) at (-7,-5) {};
				\node (v13) at (1.25,-3.75) {};
				\node (v15) at (-6.5,-6.5) {};
				\node (v16) at (-7,-2) {};
				\node (v17) at (-2,-7) {};
				\node (v18) at (-8.5,-3.5) {};
				\node (v19) at (-3.5,-8.5) {};
				\node (v14) at (-7.5,-7.5) {};
				\draw [thick] (v15) edge (v14);
				\draw [thick] (v7) edge (v10);
				\draw [line width =0.5em,white] (v5) edge (v13);
				\draw [thick] (v5) edge (v13);
				\draw [line width =0.5em,white] (v11) edge (v12);
				\draw [thick] (v11) edge (v12);
				\draw [line width =0.5em,white] (-4.75,-2.75) ellipse (1.2 and 1.2);
				\draw [line width =0.5em,white] (-2.75,-4.75) ellipse (1.2 and 1.2);
				\draw [thick] (-2.75,-4.75) node (v21) {} ellipse (1.2 and 1.2);
				\draw [thick] (-4.75,-2.75) node (v20) {} ellipse (1.2 and 1.2);
				\draw [line width =0.5em,white] (v20) edge (v10);
				\draw [thick] (v20) edge (v10);
				\draw [line width =0.5em,white] (v21) edge (v11);
				\draw [thick] (v21) edge (v11);
				\node at (-2,1.75) {$ A $};
				\node at (0.25,-0.25) {$ B $};
				\node at (-4.75,2.25) {$ T $};
				\node at (-2.75,-2) {$ I $};
				\node at (-0.5,-4) {$ J $};
				\node at (4,4) {$ S $};
				\node at (-8.5,-8.5) {$ S $};
				\draw[very thick, red]  (v1) edge (v3);
				\draw[very thick, red]  (v2) edge (v4);
				\draw[very thick]  (v1) edge (v2);
				\draw[very thick]  (v3) edge (v4);
				\node [draw,outer sep=0,inner sep=1,minimum size=10,minimum width = 26,fill=white,rotate=-45] (v9) at (1.5,1.5) {$ b $};
				\node [draw,outer sep=0,inner sep=1,minimum size=10,minimum width = 26,fill=white,rotate=-45] (v90) at (-6,-6) {$ b^* $};
				\node (v8) at (3,3) {};
				\draw [thick] (v8) edge (v9);
				\end{tikzpicture}
			\end{array}\\
		&= \sum_{\substack{A,B,T \\ b}}\frac{d(T)\Sc_{I,A}\Sc_{B,J}}{d(\C)}
			\begin{array}{c}
				\begin{tikzpicture}[scale=0.2,every node/.style={inner sep=0,outer sep=-1}]
				\node (v1) at (0.5,5.5) {};
				\node (v4) at (-2,-7) {};
				\node (v2) at (5.5,0.5) {};
				\node (v3) at (-7,-2) {};
				\node (v5) at (-4,1) {};
				\node (v6) at (1,-4) {};
				\node (v7) at (0,2) {};
				\node (v11) at (2,0) {};
				\node (v12) at (-2,-4) {};
				\node (v10) at (-4,-2) {};
				\node (v13) at (1,-4) {};
				\draw [thick] (v7) edge (v10);
				\draw [line width =0.5em,white] (v5) edge (v13);
				\draw [thick] (v5) edge (v13);
				\draw [line width =0.5em,white] (v11) edge (v12);
				\draw [thick] (v11) edge (v12);
				\node at (-2.125,1.625) {$ A $};
				\node at (0,-0.5) {$ B $};
				\node at (-5,2) {$ T $};
				\node at (4,4) {$ S $};
				\node at (-5.5,-5.5) {$ S $};
				\draw[very thick, red]  (v1) edge (v3);
				\draw[very thick, red]  (v2) edge (v4);
				\draw[very thick]  (v1) edge (v2);
				\draw[very thick]  (v3) edge (v4);
				\node [draw,outer sep=0,inner sep=1,minimum size=10,minimum width = 26,fill=white,rotate=-45] (v9) at (1.5,1.5) {$ b $};
				\node [draw,outer sep=0,inner sep=1,minimum size=10,minimum width = 26,fill=white,rotate=-45] (v90) at (-3,-3) {$ b^* $};
				\node (v8) at (3,3) {};
				\draw [thick] (v8) edge (v9);
				\node (v14) at (-4.5,-4.5) {};
				\draw [thick] (v90) edge (v14);
				\end{tikzpicture}
			\end{array} = \sum\limits_{A,B} \frac{\Sc_{I,A}\Sc_{B,J}}{d(A)d(B)} \e(S,A,B)
		\end{align*}
	where the second equality is due to a so-called `handle-slide move' (for a proof within this context see~\cite[Proposition 7.2]{hardiman2020graphical}) and the penultimate equality is due to~\cite[Lemma 3.1.4]{Bak01}.
 	\endproof

    \end{PROP}

Let $ \T $ be given by the direct sum of objects in $ \I $, i.e.\ $ \T = \bigoplus_S S $. For $ I,J \in \I $ we consider the following (diagonal) endomorphism of $ \T $,
	\begin{align*}
	\nu_I^J = \frac{1}{d(\C)} \bigoplus_S d(S) \gamma^S_{IJ}.
	\end{align*}

	\begin{THM} \label{thm:s_invariance}

	Let $ \C $ be a modular tensor category and let $ F $ be an object in $ \RTC $. Then $ F $ is $ \Sc $-invariant if and only if $ \tr F\left(\nu_I^J\right) = \tr F\left(\e_I^J\right) $ for all $ I,J \in \I $.

	\proof

	By Proposition~\ref{prop:s_matrix_appears}, we have
		\begin{align}\label{eq:apply_prop_s}
		F(\nu_I^J) = \frac{1}{d(\C)}\bigoplus_S d(S) \sum_{A,B} \frac{\Sc_{I,A}\Sc_{B,J}}{d(A)d(B)} F(\e(S,A,B)).
		\end{align}
	By the definition of $ \e(S,A,B) $, in the idempotent completion of $ \TC $ we have
		\begin{align} \label{eq:idem_isotipic_identity}
		(S,\e(S,A,B)) = S_A^B \vprod \eo{A}{B}
		\end{align}
	where $ S_A^B $ denotes the multiplicity space of the simple $ \eo{A}{B} $ in $ S $. Furthermore, by~\cite[Proposition 7.4]{hardiman2020graphical} $ S_A^B = \Hom_{\C}(S,AB) $. This implies the following sequence of natural isomorphisms,
		\begin{align} \label{eq:iso_spaces}
		F(S)_{F(\e(S,A,B))} = \Hom_{\C}(S,AB) \otimes F(AB)_{F(\e_A^B)} = \Hom_{\C}(S,AB) \otimes F_A^B
		\end{align}
 	where, for $ V $ a vector space and $ i \in \End(V) $, we write $ V_i $ to denote the corresponding subspace. As this implies $ \tr F(\e(S,A,B)) = \hom_{\C}(S,AB) \dim (F_A^B) $, we may now simply compute the trace,
		\begin{align*}
		\tr F(\nu_I^J) &= \frac{1}{d(\C)} \sum_{A,B,S} d(S) \frac{\Sc_{I,A}\Sc_{B,J}}{d(A)d(B)} \hom_{\C}(S,AB) \dim (F_A^B) \\
	 	&= \frac{1}{d(\C)} \sum_{A,B}  \left(\sum_S d(S) \hom_{\C}(S,AB)\right) \frac{\Sc_{I,A}\Sc_{B,J}}{d(A)d(B)} \dim (F_A^B) \\
		&= \frac{1}{d(\C)} \sum_{A,B} \Sc_{I,A}\Sc_{B,J} \dim (F_A^B) \\
	 	&= \left(\Sc Z(F) \Sc^{-1}\right)_{IJ}.
		\end{align*}
	As by definition $ \tr F(\e_I^J) = \dim F_I^J = Z(F)_{IJ} $, this proves the claim. \endproof

	\end{THM}

As in Section~\ref{sec:graph_char_t_inv}, we shall now rewrite this result into an equivalent formulation. Let $ X $ be in $ \TC $ and let $ F $ be in $ \RTC $. As before we consider two elements in $ \End_{\TC}(X)^* $. One is simply $ \Tr_{F,X} \colon \alpha \mapsto \tr F(\alpha) $ again, the other is given by
	\begin{align}  \label{eq:s_twisted_trace}
	\Tr_{F,X}^\Sc \colon \alpha \mapsto \sum_S \tr F \left(
	\begin{array}{c}
		\begin{tikzpicture}[scale=0.2,every node/.style={inner sep=0,outer sep=-1}]
		\node (v1) at (0,4) {};
		\node (v4) at (0,-4) {};
		\node (v2) at (4,0) {};
		\node (v3) at (-4,0) {};
		\node (v5) at (-2,2) {};
		\node (v6) at (2,-2) {};
		\node (v7) at (2,2) {};
		\node (v11) at (-2,-2) {};
		\node [draw,diamond,outer sep=0,inner sep=.5,minimum size=22,fill=white,rotate=-90] (v9) at (0,0) {$ \alpha_S $};
		\node at (3,3) {$ S $};
		\node at (-3,-3) {$ S $};
		\node at (-3,3) {$ X^\vee $};
		\node at (3,-3) {$ X^\vee $};
		\draw [thick] (v9) edge (v6);
		\draw [thick] (v9) edge (v7);
		\draw [thick] (v9) edge (v11);
		\draw [thick] (v5) edge (v9);
		\draw[very thick, red]  (v1) edge (v3);
		\draw[very thick, red]  (v2) edge (v4);
		\draw[very thick]  (v1) edge (v2);
		\draw[very thick]  (v3) edge (v4);
		\end{tikzpicture}
	\end{array}\right).
	\end{align}
where $ \alpha_S $ is the $ \Hom_{\C}(SX,XS) $ component of $ \alpha $. The following lemma gives a more intuitive description of $ \Tr_{F,X}^\Sc $.

	\begin{LEMMA} \label{lem:compute_s_twist}

	Let $ X,G $ be in $ \C $ and let $ f $ be in $ \End_{\C}(GX,XG) $. For
		\begin{align*}
		\alpha =
			\begin{array}{c}
				\begin{tikzpicture}[scale=0.2,every node/.style={inner sep=0,outer sep=-1}]
				\node (v1) at (0,4) {};
				\node (v4) at (0,-4) {};
				\node (v2) at (4,0) {};
				\node (v3) at (-4,0) {};
				\node (v5) at (-2,2) {};
				\node (v6) at (2,-2) {};
				\node (v7) at (2,2) {};
				\node (v11) at (-2,-2) {};
				\node [draw,diamond,outer sep=0,inner sep=.5,minimum size=22,fill=white,rotate=0] (v9) at (0,0) {$ f $};
				\node at (3,3) {$ X $};
				\node at (-3,-3) {$ X $};
				\node at (-3,3) {$ G $};
				\node at (3,-3) {$ G $};
				\draw [thick] (v9) edge (v6);
				\draw [thick] (v9) edge (v7);
				\draw [thick] (v9) edge (v11);
				\draw [thick] (v5) edge (v9);
				\draw[very thick, red]  (v1) edge (v3);
				\draw[very thick, red]  (v2) edge (v4);
				\draw[very thick]  (v1) edge (v2);
				\draw[very thick]  (v3) edge (v4);
				\end{tikzpicture}
			\end{array} \quad \text{we have} \quad \Stwist{\alpha} = \tr F \left(
				\begin{array}{c}
					\begin{tikzpicture}[scale=0.2,every node/.style={inner sep=0,outer sep=-1},rotate=-90]
					\node (v1) at (0,4) {};
					\node (v4) at (0,-4) {};
					\node (v2) at (4,0) {};
					\node (v3) at (-4,0) {};
					\node (v5) at (-2,2) {};
					\node (v6) at (2,-2) {};
					\node (v7) at (2,2) {};
					\node (v11) at (-2,-2) {};
					\node [draw,diamond,outer sep=0,inner sep=.5,minimum size=22,fill=white,rotate=-90] (v9) at (0,0) {$ f $};
					\node at (3,3) {$ X^\vee $};
					\node at (-3,-3) {$ X^\vee $};
					\node at (-3,3) {$ G $};
					\node at (3,-3) {$ G $};
					\draw [thick] (v9) edge (v6);
					\draw [thick] (v9) edge (v7);
					\draw [thick] (v9) edge (v11);
					\draw [thick] (v5) edge (v9);
					\draw[very thick]  (v1) edge (v3);
					\draw[very thick]  (v2) edge (v4);
					\draw[very thick, red]  (v1) edge (v2);
					\draw[very thick, red]  (v3) edge (v4);
					\end{tikzpicture}
				\end{array} \right).
		\end{align*}
	\proof

	We compute,
		\begin{align*}
		\tr F \left(
			\begin{array}{c}
				\begin{tikzpicture}[scale=0.2,every node/.style={inner sep=0,outer sep=-1},rotate=-90]
				\node (v1) at (0,4) {};
				\node (v4) at (0,-4) {};
				\node (v2) at (4,0) {};
				\node (v3) at (-4,0) {};
				\node (v5) at (-2,2) {};
				\node (v6) at (2,-2) {};
				\node (v7) at (2,2) {};
				\node (v11) at (-2,-2) {};
				\node [draw,diamond,outer sep=0,inner sep=.5,minimum size=22,fill=white,rotate=-90] (v9) at (0,0) {$ f $};
				\node at (3,3) {$ X^\vee $};
				\node at (-3,-3) {$ X^\vee $};
				\node at (-3,3) {$ G $};
				\node at (3,-3) {$ G $};
				\draw [thick] (v9) edge (v6);
				\draw [thick] (v9) edge (v7);
				\draw [thick] (v9) edge (v11);
				\draw [thick] (v5) edge (v9);
				\draw[very thick]  (v1) edge (v3);
				\draw[very thick]  (v2) edge (v4);
				\draw[very thick, red]  (v1) edge (v2);
				\draw[very thick, red]  (v3) edge (v4);
				\end{tikzpicture}
			\end{array} \right) &= \sum_{S,b} \tr F \left(
			\begin{array}{c}
				\begin{tikzpicture}[scale=0.2,every node/.style={inner sep=0,outer sep=-1},rotate=-90]
				\node (v1) at (-1.5,5.5) {};
				\node (v4) at (0,-4) {};
				\node (v2) at (4,0) {};
				\node (v3) at (-5.5,1.5) {};
				\node (v5) at (-3.5,3.5) {};
				\node (v6) at (2,-2) {};
				\node (v7) at (3,1) {};
				\node (v11) at (-1,-3) {};
				\node [draw,diamond,outer sep=0,inner sep=.5,minimum size=0,fill=white,rotate=-90] (v9) at (1,-1) {\scriptsize $ f $};
				\node at (4,2) {$ X^\vee $};
				\node at (-2,-4) {$ X^\vee $};
				\node at (-0.25,2.25) {$ S $};
				\node at (-4.5,4.5) {$ G $};
				\node at (3,-3) {$ G $};
				\draw [thick] (v9) edge (v6);
				\draw [thick] (v9) edge (v7);
				\draw [thick] (v9) edge (v11);
				\draw [thick] (v5) edge (v9);
				\node [draw,outer sep=0,inner sep=.5,minimum size=9,fill=white,rotate=-45] (v90) at (-0.5,0.5) {\scriptsize $ b $};
				\node [draw,outer sep=0,inner sep=.5,minimum size=9,fill=white,rotate=-45] (v90) at (-2.5,2.5) {\scriptsize $ b^* $};
				\draw[very thick]  (v1) edge (v3);
				\draw[very thick]  (v2) edge (v4);
				\draw[very thick, red]  (v1) edge (v2);
				\draw[very thick, red]  (v3) edge (v4);
				\end{tikzpicture}
			\end{array} \right) \\
 		&= \sum_{S,b}\tr F \left(
			\begin{array}{c}
				\begin{tikzpicture}[scale=0.2,every node/.style={inner sep=0,outer sep=-1},rotate=-90]
				\node (v1) at (-1.5,5.5) {};
				\node (v4) at (0,-4) {};
				\node (v2) at (4,0) {};
				\node (v3) at (-5.5,1.5) {};
				\node (v5) at (-3.5,3.5) {};
				\node (v6) at (2,-2) {};
				\node (v7) at (1.25,2.75) {};
				\node (v11) at (-2.75,-1.25) {};
				\node [draw,diamond,outer sep=0,inner sep=.5,minimum size=0,fill=white,rotate=-90] (v9) at (-0.75,0.75) {\scriptsize $ f $};
				\node at (2.25,3.75) {$ X^\vee $};
				\node at (-3.75,-2.25) {$ X^\vee $};
				\node at (-4.5,4.5) {$ S $};
				\node at (3,-3) {$ S $};
				\draw [thick] (v9) edge (v6);
				\draw [thick] (v9) edge (v7);
				\draw [thick] (v9) edge (v11);
				\draw [thick] (v5) edge (v9);
				\node [draw,outer sep=0,inner sep=.5,minimum size=9,fill=white,rotate=-45] (v90) at (-2.25,2.25) {\scriptsize $ b $};
				\node [draw,outer sep=0,inner sep=.5,minimum size=9,fill=white,rotate=-45] (v90) at (0.75,-0.75) {\scriptsize $ b^* $};
				\draw[very thick]  (v1) edge (v3);
				\draw[very thick]  (v2) edge (v4);
				\draw[very thick, red]  (v1) edge (v2);
				\draw[very thick, red]  (v3) edge (v4);
				\end{tikzpicture}
			\end{array} \right) = \Stwist{\alpha}
		\end{align*}
	where the second equality is due to the commutativity of $ \tr $. \endproof

	\end{LEMMA}

	\begin{THM} \label{thm:s_invariance_reformulated}

	Let $ \C $ be a modular tensor category and let $ F $ be an object in $ \RTC $. Then $ F $ is $ \Sc $-invariant if and only if $ \Tr_{F,X}^{\Sc} = \Tr_{F,X} $ for all $ X $ in $ \TC $.

	\proof 	We start start by supposing that $ \Tr_{F,X}^{\Sc} = \Tr_{F,X} $. As $ \Tr_{F,X} (\e_I^J) = \tr F(\e_I^J) $ and $ \Stwist{\e_I^J} = \tr F (\nu_I^J) $, Theorem~\ref{thm:s_invariance} allows us to conclude immediately.

	Now we suppose that $ F $ is S-invariant. We take $ X $ in $ \TC $, $ \alpha \in \End_{\TC}(X) $ and a basis $ \B \subset \Hom_{\TC}(\e_I^J,X) $. For $ b,b' \in \B $ and $ I,J \in \I $, we consider the scalar
  		\begin{align*}
  		\lambda_{IJ}^{bb'} \defeq b^* \circ \alpha \circ b' \in \End_{\TC}(\e_I^J) = \fld.
  		\end{align*}
	As $ \{\e_I^J\}_{I,J \in \I} $ forms a complete set of primitive idempotents in $ \TC $, we have
 		\begin{align*}
		\alpha = \sum_{\substack{I,J \\ b,b' }} \lambda_{IJ}^{bb'} \ b' \circ b^* = \sum_{\substack{I,J \\ b,b' }} \lambda_{IJ}^{bb'} \sum_{R,S} \begin{array}{c}
			\begin{tikzpicture}[scale=0.3,every node/.style={inner sep=0,outer sep=-1}]
			\node (v1) at (-1.375,2.875) {};
			\node (v4) at (-3.375,-3.125) {};
			\node (v2) at (0.625,0.875) {};
			\node (v3) at (-5.375,-1.125) {};
			\node (v9) at (-0.375,1.875) {};
			\node (v6) at (-4.375,-2.125) {};
			\node at (0.25,2.5) {$ X $};
			\node at (-5,-2.75) {$ X $};
			\node at (-5,0.75) {$ R $};
			\node at (-3.25,2.5) {$ S $};
			\node at (-1.5,-2.75) {$ R $};
			\node at (0.25,-1) {$ S $};
			\node (v11) at (-3.5,-0.75) {};
			\node (v12) at (-3,-1.25) {};
			\draw [thick] (v11) edge (v12);
			\node (v14) at (-1.75,1) {};
			\node (v13) at (-1.25,0.5) {};
			\draw [thick] (v13) edge (v14);
			\node (v15) at (-4.25,0) {};
			\node (v16) at (-2.5,1.75) {};
			\draw [thick] (v15) to (v11);
			\draw [thick] (v16) to (v14);
			\node (v17) at (-0.5,-0.25) {};
			\node (v18) at (-2.25,-2) {};
			\draw [thick] (v13) to (v17);
			\draw [thick] (v12) to (v18);
			\draw [thick] (v14) edge (v11);
			\draw [thick] (v13) edge (v12);
			\node [rotate=-45,draw,outer sep=0,inner sep=0,minimum size=10,fill=white] (v50) at (-3.25,-1) {$\scriptstyle b'_R $};
			\node [rotate=-45,draw,outer sep=0,inner sep=0,minimum size=10,fill=white] (v5) at (-1.5,0.75) {$\scriptstyle b_S^* $};
			\draw [thick] (v50) edge (v6);
			\draw[thick]  (v9) edge (v5);
			\draw[very thick,red]  (v1) edge (v3);
			\draw[very thick,red]  (v2) edge (v4);
			\draw[very thick]  (v1) edge (v2);
			\draw[very thick]  (v3) edge (v4);
			\end{tikzpicture}
		\end{array}.
 		\end{align*}
	and
		\begin{align*}
		\Tr_{F,X} (\alpha) = \sum_{\substack{I,J \\ b,b' }} \lambda_{IJ}^{bb'} \tr F (b' \circ b^*) = \sum_{\substack{I,J \\ b,b' }} \lambda_{IJ}^{bb'} \tr F (b^* \circ b') = \sum_{\substack{I,J \\ b}} \lambda_{IJ}^{bb} \tr F (\e_I^J).
		\end{align*}
	By Lemma~\ref{lem:compute_s_twist}, we can compute $ \Stwist{\alpha} $ as
		\begin{align*}
		& \hspace{1.4em} \sum_{\substack{I,J \\ b,b' }} \lambda_{IJ}^{bb'} \sum_{R,S} \tr F \left(
			\begin{array}{c}
				\begin{tikzpicture}[scale=0.3,every node/.style={inner sep=0,outer sep=-1},rotate=-90]
				\node (v1) at (-1.375,2.875) {};
				\node (v4) at (-3.375,-3.125) {};
				\node (v2) at (0.625,0.875) {};
				\node (v3) at (-5.375,-1.125) {};
				\node (v9) at (-0.375,1.875) {};
				\node (v6) at (-4.375,-2.125) {};
				\node at (0.25,2.5) {$ X^\vee $};
				\node at (-5,-2.75) {$ X^\vee $};
				\node at (-5,0.75) {$ R $};
				\node at (-3.25,2.5) {$ S $};
				\node at (-1.5,-2.75) {$ R $};
				\node at (0.25,-1) {$ S $};
				\node (v11) at (-3.5,-0.75) {};
				\node (v12) at (-3,-1.25) {};
				\draw [thick] (v11) edge (v12);
				\node (v14) at (-1.75,1) {};
				\node (v13) at (-1.25,0.5) {};
				\draw [thick] (v13) edge (v14);
				\node (v15) at (-4.25,0) {};
				\node (v16) at (-2.5,1.75) {};
				\draw [thick] (v15) to (v11);
				\draw [thick] (v16) to (v14);
				\node (v17) at (-0.5,-0.25) {};
				\node (v18) at (-2.25,-2) {};
				\draw [thick] (v13) to (v17);
				\draw [thick] (v12) to (v18);
				\draw [thick] (v14) edge (v11);
				\draw [thick] (v13) edge (v12);
				\node [rotate=-135,draw,outer sep=0,inner sep=0,minimum size=10,fill=white] (v50) at (-3.25,-1) {$\scriptstyle b'_R $};
				\node [rotate=-135,draw,outer sep=0,inner sep=0,minimum size=10,fill=white] (v5) at (-1.5,0.75) {$\scriptstyle b^*_S$};
				\draw [thick] (v50) edge (v6);
				\draw[thick]  (v9) edge (v5);
				\draw[very thick]  (v1) edge (v3);
				\draw[very thick]  (v2) edge (v4);
				\draw[very thick, red]  (v1) edge (v2);
				\draw[very thick, red]  (v3) edge (v4);
				\end{tikzpicture}
			\end{array} \right)
		=\sum_{\substack{I,J \\ b,b' }} \lambda_{IJ}^{bb'}  \sum_{R,S} \tr F \left(
			\begin{array}{c}
				\begin{tikzpicture}[scale=0.3,every node/.style={inner sep=0,outer sep=-1},rotate=-90]
				\node (v1) at (-2,3.5) {};
				\node (v4) at (-2.75,-3.75) {};
				\node (v2) at (1.25,0.25) {};
				\node (v3) at (-6,-0.5) {};
				\node (v9) at (-0.4375,1.8125) {};
				\node (v6) at (-4.3125,-2.0625) {};
				\node at (-5.25,-2.75) {$ X^\vee $};
				\node at (-4,-4) {$ R^\vee $};
				\node at (-6.5,-1.5) {$ R $};
				\node at (0.25,2.75) {$ X^\vee $};
				\node at (1.5,1.5) {$ R^\vee $};
				\node at (-1,4) {$ R $};
				\node at (-5.625,1.375) {$ S $};
				\node at (-3.875,3.125) {$ R $};
				\node at (-0.875,-3.375) {$ S $};
				\node at (0.875,-1.625) {$ R $};
				\node (v11) at (-3.5,-0.75) {};
				\node (v12) at (-3,-1.25) {};
				\draw [thick] (v11) edge (v12);
				\node (v14) at (-1.75,1) {};
				\node (v13) at (-1.25,0.5) {};
				\draw [thick] (v13) edge (v14);
				\node (v15) at (-5.4375,-0.9375) {};
				\node (v16) at (-4.6875,0.6875) {};
				\draw [thick] (v15) to[out=55,in=135] (v11);
				\draw [thick] (v16) to[out=-45,in=135] (v14);
				\node (v17) at (-1.8125,-2.6875) {};
				\node (v18) at (-3.1875,-3.1875) {};
				\draw [thick] (v13) to[out=-45,in=135] (v17);
				\draw [thick] (v12) to[in=50,out=20] (v18);
				\draw [thick] (v14) edge (v11);
				\draw [thick] (v13) edge (v12);
				\node [rotate=-135,draw,outer sep=0,inner sep=0,minimum size=10,fill=white] (v50) at (-3.25,-1) {$\scriptstyle b'_R $};
				\node [rotate=-135,draw,outer sep=0,inner sep=0,minimum size=10,fill=white] (v5) at (-1.5,0.75) {$\scriptstyle b^*_S$};
				\draw [thick] (v50) edge (v6);
				\draw[thick]  (v9) edge (v5);
				\node (v7) at (-1.3125,2.6875) {};
				\node (v8) at (-3.25,2.1667) {};
				\draw [thick] (v7) to[out=-130,in=-45] (v8);
				\node (v10) at (0.4375,0.9375) {};
				\node (v19) at (-0.0625,-0.9375) {};
				\draw [thick] (v10) to[out=-130,in=135] (v19);
				\draw[very thick]  (v1) edge (v3);
				\draw[very thick]  (v2) edge (v4);
				\draw[very thick, red]  (v1) edge (v2);
				\draw[very thick, red]  (v3) edge (v4);
				\end{tikzpicture}
			\end{array}
		\right) \\
		&= \sum_{\substack{I,J \\ b,b' }}\lambda_{IJ}^{bb'} \sum_{R,S} \tr F \left(  \begin{array}{c}
			\begin{tikzpicture}[scale=0.3,every node/.style={inner sep=0,outer sep=-1},rotate=-90]
			\node (v1) at (-1.5,2.75) {};
			\node (v4) at (-3.25,-3) {};
			\node (v2) at (0.5,0.75) {};
			\node (v3) at (-5.25,-1) {};
			\node (v9) at (-0.5,1.75) {};
			\node at (-0.5,3) {$ I^\vee $};
			\node at (0.75,1.75) {$ J^\vee $};
			\node at (-5.5,-2) {$ I^\vee $};
			\node at (-4.25,-3.25) {$ J^\vee $};
			\node at (-5,0.75) {$ S $};
			\node at (-3.25,2.5) {$ R $};
			\node at (-1.5,-2.75) {$ S $};
			\node at (0.25,-1) {$ R $};
			\node (v11) at (-3.5,-0.75) {};
			\node (v12) at (-3,-1.25) {};
			\draw [thick] (v11) edge (v12);
			\node (v14) at (-1.75,1) {};
			\node (v13) at (-1.25,0.5) {};
			\draw [thick] (v13) edge (v14);
			\node (v15) at (-4.25,0) {};
			\node (v16) at (-2.5,1.75) {};
			\draw [thick] (v15) to (v11);
			\draw [thick] (v16) to (v14);
			\node (v17) at (-0.5,-0.25) {};
			\node (v18) at (-2.25,-2) {};
			\draw [thick] (v13) to (v17);
			\draw [thick] (v12) to (v18);
			\node (v8) at (-4,-2.25) {};
			\node (v7) at (-4.5,-1.75) {};
			\node (v6) at (-0.75,2) {};
			\node (v10) at (-0.25,1.5) {};
			\draw [thick] (v11) edge (v7);
			\draw [thick] (v12) edge (v8);
			\draw [thick] (v14) edge (v6);
			\draw [thick] (v13) edge (v10);
			\node [rotate=-135,draw,outer sep=0,inner sep=0,minimum size=10,fill=white] (v50) at (-3.25,-1) {$\scriptstyle b^*_S $};
			\node [rotate=-135,draw,outer sep=0,inner sep=0,minimum size=10,fill=white] (v5) at (-1.5,0.75) {$\scriptstyle b'_R $};
			\draw [thick] (v50) edge (v5);
			\draw[very thick]  (v1) edge (v3);
			\draw[very thick]  (v2) edge (v4);
			\draw[very thick,red]  (v1) edge (v2);
			\draw[very thick,red]  (v3) edge (v4);
			\end{tikzpicture}
		\end{array} \right) = \sum_{\substack{I,J \\ b }}\lambda_{IJ}^{bb} \sum_{T} \frac{d(T)}{d(\C)} \tr F \left( \begin{array}{c}
			\begin{tikzpicture}[scale=0.15,every node/.style={inner sep=0,outer sep=-1},rotate=-90]
			\node (v1) at (0,5) {};
			\node (v4) at (0,-5) {};
			\node (v2) at (5,0) {};
			\node (v3) at (-5,0) {};
			\node (v5) at (-2.5,2.5) {};
			\node (v6) at (2.5,-2.5) {};
			\node (v7) at (1.5,3.5) {};
			\node (v11) at (3.5,1.5) {};
			\node (v12) at (-1.5,-3.5) {};
			\node (v10) at (-3.5,-1.5) {};
			\node (v13) at (2.5,-2.5) {};
			\draw [thick] (v7) edge (v10);
			\draw [line width =0.5em,white] (v5) edge (v13);
			\draw [thick] (v5) edge (v13);
			\draw [line width =0.5em,white] (v11) edge (v12);
			\draw [thick] (v11) edge (v12);
			\node at (-5.25,-2.5) {$ I^\vee $};
			\node at (-2.75,-5.25) {$ J^\vee $};
			\node at (-4,4) {$ T $};
			\node at (4,-4) {$ T $};
			\draw[very thick]  (v1) edge (v3);
			\draw[very thick]  (v2) edge (v4);
			\draw[very thick, red]  (v1) edge (v2);
			\draw[very thick, red]  (v3) edge (v4);
			\end{tikzpicture}
		\end{array} \right)\\
		&= \sum_{\substack{I,J \\ b}} \lambda_{IJ}^{bb} \tr F (\nu_I^J) =  \sum_{\substack{I,J \\ b}} \lambda_{IJ}^{bb} \tr F (\e_I^J) = \Tr_{F,X}(\alpha)
		\end{align*}
	where the first equality is due to the invariance of $ \tr $ under conjugation by an isomorphism, the third is due to the fact that $ b^* \circ b' = \delta_{b,b'} \e_I^J $ and the fifth is due to Theorem~\ref{thm:s_invariance}. \endproof

	\end{THM}

	\section{Conclusion} \label{sec:conclusion}

We are now in position to combine the results of Section~\ref{sec:graph_char_t_inv} and Section~\ref{sec:graph_char_s_inv} into our main result.

	\begin{COR} \label{cor:main_result}

	Let $ F $ be in $ \RTC $. Then $ F $ is modular invariant if and only if
		\begin{align*}
		\Tr_{F,X} = \Tr_{F,X}^{\Tc} = \Tr_{F,X}^{\Sc}
		\end{align*}
	for all $ X $ in $ \TC $.

	\proof This follows immediately from Corollary~\ref{cor:t_invariance_reformed} and Theorem~\ref{thm:s_invariance_reformulated}. \endproof

	\end{COR}

Taken together,\! $ \Tr_{F,X}^{\Tc} $ and $ \Tr_{F,X}^{\Sc} $ are extremely telling. In particular, their image under $ \alpha \in \End_{TC}(X) $  is given by taking the standard trace $ \Tr_{F,X} $ and pre-composing it with an action of the standard generators of the mapping class group of the torus on the graphical description of $ \alpha $\footnote{This is clear for $ \Tr_{F,X}^{\Tc} $ and for $ \Tr_{F,X}^{\Sc} $ it is the content of Lemma~\ref{lem:compute_s_twist}.}. Therefore Corollary~\ref{cor:main_result} may be reformulated into stating that modular invariance of an object in $ \RTC $ amounts to a form of invariance under the mapping class group of the torus, giving us a characterization which doesn't distinguish certain generators.

Furthermore, this picture suggests a deeper relationship with the motivational mathematical physics described in the introduction. We recall that the modular invariance of the partition function encodes the conformal invariance of the field theory. We also recall that the partition function gives the value of the theory on a torus parametrised (in the standard way) by a point in the complex upper half-plane $ \mathbb{H} $. As $ \mathbb{H} $ is the Teichm\"uller space, taking the quotient by the mapping class group gives the moduli space. Therefore conformal invariance with respect to the torus also comes down to a form of invariance under the mapping class group of the torus. Drawing out this correspondance and connecting these two actions of the mapping class group is the subject of ongoing reseach.

\bibliographystyle{alpha}
\bibliography{graphical_s_inv}

\end{document}